\documentstyle[12pt,bezier]{article}

\newcommand{\sect}[1]{\setcounter{equation}{0}\section{#1}}

\newcommand{\eq}{\begin{equation}}
\newcommand{\eqa}{\begin{eqnarray}}
\newcommand{\en}{\end{equation}}
\newcommand{\ena}{\end{eqnarray}}
\newcommand{\enn}{\nonnumber \end{equation}}
\newcommand{\spz}{\hspace{0.7cm}}
\newcommand{\ie}{{\it i.e.\ }}
\def\err{\ell^{\cal R}}

\def\sk{\vskip .4cm}
\def\noi{\noindent}
\def\om{\omega}
\def\al{\alpha}

\def\Ga{\Gamma}

\def\del{\delta}

\def\epsi{\varepsilon}

\def\part{\partial}

\def\f#1#2{ f_{#2}{}^{#1} }

\def\M#1#2{ M_{#1}^{~#2} }

\def\D{\Delta}

\def\invG{{}_{inv}\Ga}
\def\invX{{}_{inv}\Xi}

\def\dia{{\scriptscriptstyle \Box}}
\def\cvd{\rightline{$\Box\!\Box\!\Box$}\sk}
\def\FF#1#2{O_{#1}{}^{#2} }

\def\DL{\Delta_{\scriptstyle \Ga}}

\def\DR{{}_{\scriptstyle \Ga}\Delta}

\def\DS{\Delta_{\scriptstyle \Xi}}

\def\DD{{}_{\scriptstyle \Xi}\Delta}

\def\N#1#2{N^{#1}{}_{#2}}
\def\vart{\vartheta}

\begin{document}
\begin{titlepage}
\rightline{IFUP-TH 15/95}
\rightline{LMU-TPW 94-14}
\rightline{q-alg/9505023}
\sk
\sk
\begin{center}{\bf VECTOR FIELDS ON QUANTUM GROUPS }\\[4em]
Paolo Aschieri${}^1$ and Peter Schupp${}^2$\\[2em]
{\sl ${}^1$Scuola Normale Superiore \\
Piazza dei Cavalieri 7, 56100 Pisa, Italy\\
and\\
Istituto Nazionale di Fisica Nucleare, sezione di Pisa, Italy } \\[2em]
{\sl ${}^2$Sektion Physik der Ludwig-Maximilians-Universit\"at M\"unchen\\
Theoretische Physik --- Lehrstuhl Professor Wess\\
Theresienstr. 37, D-80333 M\"unchen } \\[6em]
\end{center}
\begin{abstract}
We construct the space of vector fields on a generic quantum group.
Its elements are products of elements of the quantum group
itself with left invariant vector fields.
We study the duality between vector fields and 1-forms and generalize
the construction to tensor fields.
A Lie derivative along any (also non left invariant) vector field
is proposed and a puzzling ambiguity in its definition discussed.
These results hold for a generic Hopf algebra.
\end{abstract}

\vskip 3cm

\noi \hrule
\vskip.2cm
\hbox{{\small{\it e-mail: }}{\small Aschieri@ux2sns.sns.it,
Schupp@lswes8.ls-wess.physik.uni-muenchen.de}}
\end{titlepage}
\newpage
\setcounter{page}{1}

\sect{Introduction}

Following the program of generalizing the differential geometry
structures to the non-commutative case, we construct on a Hopf algebra the
analog of the space of vector fields.

Indeed in the literature the quantum Lie algebra of left invariant vector
fields as well as the space of 1-forms
has been extensively analyzed \cite{Wor,Jurco,Zumino,
Watamura,Su,SWZ2,AC},
while the notion of generic vector field on a Hopf
algebra and the duality relation with the space of $1$-forms deserves
more study \cite{SWZ3,SW}.

We will see how left invariant vector fields
generate the whole space of vector fields. This space can be
also characterized as the bicovariant bimodule (vector bundle)
dual to that of 1-forms.

Throughout this paper we will deal  with a Hopf algebra $A$
\cite{Swe}
over $\mbox{\boldmath$C$}$ with coproduct $\D\::~A\rightarrow A\otimes A
$, counit $\epsi\::~A\rightarrow \mbox{\boldmath$C$}$ and invertible
antipode
$S\::~A \rightarrow A  $.\\
Particular cases of Hopf algebras are quantum  groups, which for
us will be Hopf algebras with one (or more) continuous parameter $q$ ;
when
$q$=$1$ the product ``${\displaystyle \cdot} $'' in $A$ becomes
commutative and
we obtain the algebra of functions on a group{}. When we will speak of
 commutative case we will refer to the Hopf algebra
$C^{\infty}(G)$ of smooth functions on a (compact) Lie group $G$.

In  Section 2 we briefly recall how to associate a Connes-type
differential calculus to a given Hopf algebra and we emphasize the role
played by the tangent vectors.\\
This construction will be effected along the lines of Woronowicz' work
\cite{Wor}.
Indeed the results in \cite{Wor} apply also to a general Hopf
algebra with invertible antipode (not necessary a compact matrix
pseudogroup). This can be shown by checking that all the formulae
used for the construction --- collected in the
appendix  of \cite{Wor} ---
hold also in the case of a Hopf algebra with invertible antipode
\footnote{Formula (A.22) in \cite{Wor} is the most difficult to prove
and necessitates the further axiom of the invertibility of the
antipode $S$; also the invertibility of the map $s$ in (A.18) relies on
the existence of $S^{-1}$. All the other formulae are direct
consequences of the Hopf algebra axioms.}.

While following the work of \cite{Wor} in spirit, we however
decided on a small but important change of conventions \cite{Zumino}
in this presentation.
It is based on the following observation: while not a priori obvious, it is
indeed possible to write all defining relations of the differential
calculus as {\em deformed commutation relations}. Now, given that this is
possible, we would like to have operators acting from the left
and to the right as one is used to. This made the change in conventions
necessary. Convertation to the old conventions (denoted by a ``W''
subscript) is possible in two
ways: either one substitutes
$$\chi_i = -S^{-1}({{}_{{}_{W\!}}}\chi_i)
, ~~ f_j{}^i = S^{-1}({{}_{{}_{W\!}}}f^i{}_j),
\spz {\it etc.}$$
into our equations as will be mentioned in the text,
or one simply reads all equations
``from the left to the right''. Notice that the ad-invariant right ideal $R$,
the differential and the forms are the same.
For quick reference we would like
to point out the first version \cite{Paolo} of this article, where
the old convention was used.

In Section 3 we construct
the space of vector fields, while in Section 4
we study the
action of the Hopf algebra on the vector fields; i.e.
we will study the push-forward of vector fields on Hopf algebras.
Then we deal with (covariant and contravariant) tensor fields
and wedge products.

In the last three sections we introduce and discuss
Lie derivatives and a contraction operator on differential forms
along generic vector fields.
These two operators are basic tools for the formulation of
deformed gravity theories \cite{Castellani},
where the relevant Lie algebra is the
$q$-Poincar\'{e} Lie algebra.

\sect{Differential Geometry on Hopf Algebras}

In the commutative case, given the
differential calculus on a (compact) Lie
group $G$, we can consider the subspace in the space of all
smooth functions $f~:~~G\longrightarrow \mbox{\boldmath$C$}$ defined by:
\eq
R\equiv\{h\in C^{\infty}(G)~ |~~~ h(1_G)=0 \mbox{ and } dh(1_G)=0 \}~,
\label{ideale}
\en
where $1_G $ is the unit of the group.\\
$R$ is a particular ideal of the Hopf algebra $C^{\infty}(G)$ namely
$({\rm ker}\epsilon)^2$.
All the information about the differential calculus  on $G$ is  contained in
$R$.\\
Indeed the space of tangent vectors at the origin of the group is
given by all the linear functionals which annihilate $R$ and any constant
function.
Locally we write a basis as $\{\partial_i|_{1_G}\}$.
Once we have this basis,
using the tangent map (namely
$TL_{g}$) induced by the left multiplication of the group on itself:
$L_gg'=gg'\;,~\,\forall g,g' \in G$ we can construct a basis of left
invariant vector fields $\{t_i\}$. The action of these vector fields
on a generic function $a$ on the group manifold is (here and in
the following: $\Delta(a) \equiv a_1 \otimes a_2$ is the coproduct
in Sweedler's notation \cite{Swe})
$$
t_i(a) = a_1 (\partial_i a_2 |_{1_G}) \equiv \partial_i|_{1_G} * a
$$
in compliance with the following picture:\\
\unitlength=1.00mm
\special{em:linewidth 0.4pt}
\linethickness{0.4pt}
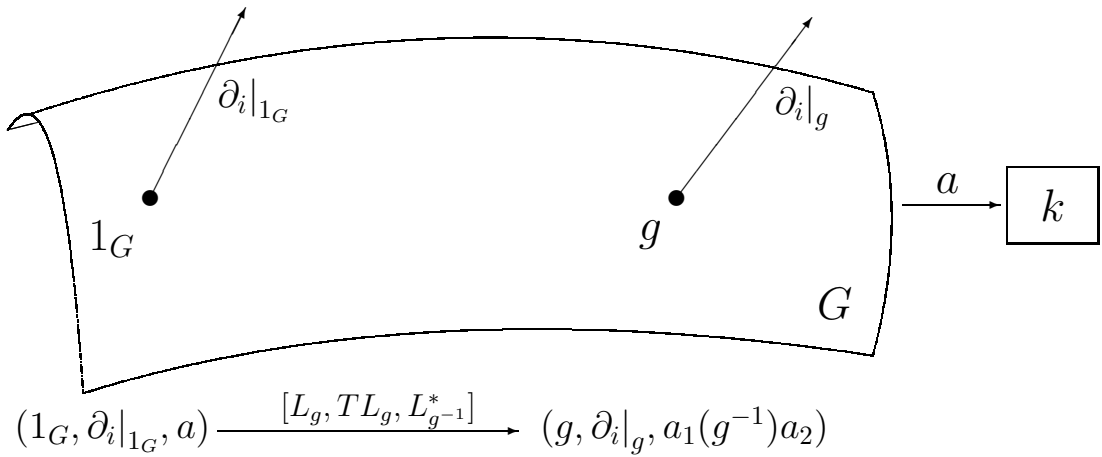
\begin{picture}(144.00,68.00)
\put(18.00,38.00){\circle*{2.00}}
\put(88.00,38.00){\circle*{2.00}}
\put(18.00,38.00){\vector(1,2){12.67}}
\put(88.00,38.00){\vector(3,4){18.00}}
\put(109.00,24.00){\makebox(0,0)[cc]{\Large$G$}}
\put(16.00,35.00){\makebox(0,0)[rt]{\Large$1_G$}}
\put(86.00,35.00){\makebox(0,0)[rt]{\Large$g$}}
\put(27.00,51.00){\makebox(0,0)[lc]{\large${\partial_i|}_{1_G}$}}
\put(101.00,50.00){\makebox(0,0)[lc]{\large${\partial_i|}_g$}}
\put(89.00,7.00){\makebox(0,0)[cc]{\large$(g,{\partial_i|}_g,a_1(g^{-1})a_2)$}}
\put(26.00,7.00){\makebox(0,0)[rc]{\large$(1_G,{\partial_i|}_{1_G},a)$}}
\put(49.00,9.00){\makebox(0,0)[cb]{[$L_g,TL_{g},L^*_{g^{-1}}$] }}
\put(132.00,32.00){\framebox(12.00,10.00)[cc]{\Large$k$}}
\put(124.00,40.00){\makebox(0,0)[cc]{\Large$a$}}
\bezier{144}(114.00,52.00)(119.00,36.00)(114.00,17.00)
\bezier{236}(-1.00,47.00)(6.00,58.00)(9.00,12.00)
\bezier{432}(9.00,12.00)(54.00,26.00)(114.00,17.00)
\bezier{472}(1.50,49.00)(59.00,68.00)(114.00,52.00)
\put(27.00,7.00){\vector(1,0){40.00}}
\put(-1.00,47.00){\line(4,1){4.07}}
\put(118.00,37.00){\vector(1,0){13.00}}
\end{picture}\\
Note that $L^*_g(a)(h) := a(gh) = a_1(g) a_2(h)$.

A generic $1$-form can be written $\rho=\om^i f_i~~[f_i\in
C^{\infty}(G)]$ where $\{\om^i\}$ is the dual basis of
$\{t_i\}$.
Finally, the differential on functions is
\eq
d=\om^i t_i ~~~~\mbox{ that is } ~~~~ df=\om^i t_i(f)~.\label{donf}
\en
\sk
In \cite{Wor} the quantum analog of $R\subset A$ is studied.
Similarly to the classical case $R$ is a right ideal of $A$,
formally $R A \subset R$, it is Ad-invariant, \ie
$\forall r \in R: r_2 \otimes S(r_1) r_3 \in R \otimes A$,
and its elements have vanishing counit.
Given a Hopf algebra $A$, it turns out that we can always find an $R$ and
construct a differential calculus (in general not  unique).
\sk
The space of tangent vectors on $A$ is then defined as:%
\footnote{There are two tangent spaces:\\ The one we have chosen
and the one corresponding to $\chi(R) = 0$ \cite{Wor}.}
\[
T\equiv\{\chi\::~A\rightarrow \mbox{\boldmath$C$}~\:|~~~ \chi ~ \mbox{ linear
functionals, }~\chi (I)=0~
\mbox{ and }~ \chi(S R)=0 \} ~,
\]
where $I$ is the unit of $A$ (in the commutative case it is the constant
function \mbox{$I(g)=1$}{} ${}\forall g\in G$). In the sequel $\{I\}$ is the
linear span of $I$.
Let $\{\chi_i\}~i=1,\,\ldots,n$ be a basis of $T$.
Consider the linear space $X$ such that
\eq
A=X\oplus R\oplus \{I\} \label{axri}
\en
$X$ is maximal in the (ordered) set of all linear subspaces of $A$ disjoint
from $R\oplus\{I\}$. From (\ref{axri}) it follows  that the dual vector space
$X^*$ is isomorphic to $T$ and therefore there are $n$ elements
$x^i\in X\subset {\mbox{ker}}\epsi$ uniquely defined by the duality
\eq
\chi_i(x^j)=\delta^j_i ~.\label{chix}
\en
Note that $\epsi(x^i)=0$ since
 $X\subset\mbox{ker}\epsi$ because
$A={\mbox{ker}}\epsi \oplus\{I\}.$
%
%
We can then define the $n^2$ linear functionals\footnote{Relation
to the conventions of \cite{Wor} (denoted by a ``W'' subscript):\\ $\chi_i =
-S^{-1}{{}_{{}_{W\!}}}\chi_i,~~
\f{i}{j} = S^{-1}{{}_{{}_{W\!}}}f{}^i{}_j,~~R
= {{}_{{}_{W\!}}}R,~~ x^i = -S{{}_{{}_{W\!}}}x^i.$}
$\f{i}{j}~:~~A\longrightarrow \mbox{\boldmath$C$}$
\eq
\!\!\!\!\!\!\!\!\!\!\!\!\!\!\!\!\! \forall a \in
A~~~~~~~~\f{i}{j}(a)\equiv\chi_j(ax^i)~,
\en
and (\ref{chix})] and we have (from
$A S(R) = S(R)$)
\eq
\chi_i(ab)=\chi_i(a)\epsi(b) + \f{j}{i}(a)\chi_j(b)~.
\label{22bis}
\en
This is the deformed Leibniz rule for the operators $\chi_i$.
In the $q=1$ case, when $R$ becomes the set defined in (\ref{ideale}),
we have $\chi_i=\partial_i|_{1_G} $ ,
 $\f{i}{j}=\delta^i_j\epsi$ and we write (\ref{22bis}) as $
\partial_i(fh)|_{1_G}=(\partial_if|_{1_G})h(1_G) + f(1_G)(\partial_i
h|_{1_G})$.
\sk
For consistency with (\ref{22bis})
the $\f{i}{j}$ must satisfy the conditions:
\begin{eqnarray}
& & \f{i}{j} (ab)= \f{k}{j} (a) \f{i}{k} (b) \label{propf1}\\
& & \f{i}{j} (I) = \del^i_j ~.\label{propf2}
\end{eqnarray}
The space of left invariant  vector fields $\invX$ is easily
constructed from $T$. Using the coproduct $\D$ we define ---
as in the commutative case ---
$\chi*a=(id\otimes \chi)\D(a) \equiv a_1 \chi(a_2)$ and
\eq
{}_{inv}\Xi \equiv \{t~ |~~ t=\chi * ~ \mbox{ where } \chi \in T\}
\label{leftXi}
\en
\noi There is a one to one correspondence $\chi_i \leftrightarrow t_i=
\chi_i*$. In order to obtain $\chi_i$ from $\chi_i*$ we simply apply
$\epsi$:
\eq
(\epsi \circ t_i)(a) = \epsi(a_1)\chi_i(a_2) = \chi_i(a),
\en
where we have used the Hopf algebra axioms.

\noi ${}_{inv}\Xi$ is the vector subspace of all linear maps from $A$ to $A$
that
is isomorphic to $T$.


We have chosen this perspective to introduce the space of left invariant
vector fields in order to point out that (also in the case of a
general Hopf algebra) it has an existence  on its own,
independent of  the space of 1-forms.
\sk
The space of $1$-forms $\Ga$ is formed by all the elements $\rho$ that
are
 written as formal products and sums of the type
\begin{eqnarray}
\rho= \om^i a_i~. \label{rhoaom}
\end{eqnarray}
\noi Here $a_i \in A$ and {$\om^i\,$} $i=1,\,\ldots,n$ is the basis dual to
$\{t_i\}$ .
We express this duality with a bracket:
\eq
\langle\chi_j,\om^i\rangle = \delta^i_j \label{braket} ~~.
\en
Relation (\ref{rhoaom}) tells us that the space of $1$-forms is freely
generated by the elements $\om^i$. By definition any $\rho$ is
decomposable in a unique way as $\rho=\om^ia_i$ and $\Ga$ is a right
$A$-module with the trivial product $(\om^ia_i)b\equiv \om^i(a_ib)$.
$\Ga$ is also a left $A$-module with the following left product:
\begin{eqnarray}
\!\!\!\!\!\!\!\!\!\!\!\!\!\!\!\!\!\!\!
\forall b\in A &~~~~~~~ & b\om^i = \om^j (\f{i}{j} * b)  \equiv \om^j
(id \otimes \f{i}{j}) \Delta (b)~. \label{omb}
\end{eqnarray}
{}From this relation it follows that:
\begin{eqnarray}
\!\!\!\!\!\!\!\!\!\!\!\!\!\!\!\!\!\!\!
\forall a\in A&~~~~~~~ & \om^ia= [(\f{i}{j} \circ S)* a]\om^j
= a_1 \f{i}{j}(Sa_2) \om^j \label{aom}
\end{eqnarray}
and that any $\rho$ can be written in a unique way in the form
\eq
\rho=b_i\om^i \label{rhoomb}
\en
 with $b_i\in A$.
\sk
Finally, the differential operator $~d~:~~A\longrightarrow \Ga~$ can be
defined through the relation:
\eq
\!\!\!\!\!\!\!\!\!\!\!\!\!\!\!\!\!\!\!\!\!\!\!\!\!\!\!\!\!\!\!\!\!
\forall a\in A ~~~~~~~~~~da = \om^i(\chi_i * a) ~.\label{dachi}
\en
Notice that this can be rewritten:
\eq
da = \om^i(\chi_i * a) = S f_i{}^j * (\chi_j * a) \om^i =
(-S \chi_i * a) \om^i
\en
where we have
used (\ref{aom}) and (\ref{22bis}), i.e. ${{}_{{}_{W\!}}}\chi_i = -
S\chi_i$.
As a consequence, the dif\-ferential calculus obtained has the following
properties \cite{Wor}:\hfill \\
\sk
\noi i) The differential operator satisfies the Leibniz rule
\eq
d(ab)=(da)b+a(db) ~~\forall a,b\in A . \label{Leibniz}
\en
Moreover  any $\rho \in \Ga$ can be expressed as
\eq
\rho=da_{\al} b_{\al} \label{adb}
\en
\noi for some $a_{\al},b_{\al}$ belonging to $A$. (Use $\om^i = d x_2^i
S^{-1} x_1^i$ below).
\sk
\noi {\sl Remark}. Once we know the operator $d$, the space of tangent
vectors on $A$, like in the commutative case, can be defined as:
\eq
T=\{\chi~ |~ \chi(a)=0 \mbox{ if and only if } Pda=0\}
\en
where $Pda \equiv da_2 S^{-1}a_1$ with
$\Delta(a)=a_{1}\otimes a_{2}.$
The linear map $P$
is a projection operator; to a given
form $\rho = \om^ia_i$ it associates
the form $P(\rho) = \epsi (a_i)\om^i$ and in particular
$P(d x^i) = \om^i$.
In the commutative case  $\epsi(a_i)$ is the value
that $a_i \in A=C^{\infty}(G)$ takes in the origin $1_G$ of the
Lie group $G$. $P(\rho)$ is then the left invariant 1-form
whose  value in the origin $1_G$ of the Lie group
equals the value of the 1-form $\rho$ in $1_G$.\footnote{In \cite{Wor}:
${{}_{{}_{W\!}}}P d a
= S(a_1) d a_2$,
such that $\om^i = P d x^i = {{}_{{}_{W\!}}}P d {{}_{{}_{W\!}}}x^i
= {{}_{{}_{W\!}}}\om^i$.} \hfill \\
\sk
\noi ii) The differential calculus is called bicovariant
because using $d$ and
the coproduct $\D$ we can define two linear compatible maps $\DL$ and $\DR$
\begin{eqnarray}
&&\DL(da)=(id \otimes d)\D(a),~~~\DL:\Ga \rightarrow A \otimes \Ga
{}~~~{\rm (left~covariance)}
\label{leftco}\\
&&\DR(da)=(d \otimes id)\D(a),~~~\DR:\Ga \rightarrow \Ga \otimes A
{}~~~{\rm (right~covariance)}
\label{rightco}\\
&&\DL(a\rho b)=\D (a)\DL(\rho)\D(b), ~~~ \DR(a\rho b)=\D(a)\DR(\rho)\D(b)
\label{Dprop0}
\end{eqnarray}
which represent the left and right action of the Hopf algebra $A$ on $\Ga$.
In the commutative case they express the
pull-back on 1-forms induced by the left or right multiplication of the
group on itself \cite{AC}.
$\DL$ and $\DR$ are compatible in the sense that
$(id \otimes \DR) \DL = (\DL \otimes id) \DR $.
In the commutative case  this formula tells us that the left and right
actions of the group on $\Ga$ commute:
$R^*_g L^*_{g'}=L^*_{g'}R^*_g ~~\forall \:g,g' \in G.$
{}From the definitions (\ref{leftco}) and (\ref{rightco}) one deduces the
 following properties \cite{Wor}:
\eqa
(\epsi \otimes id) \DL (\rho)=\rho,& & ~~~ (id \otimes \epsi) \DR (\rho)=
\rho \label{Dprop1}\\
(\D \otimes id)\DL=(id\otimes\DL)\DL,& & ~~~ (id\otimes\D)\DR=(\DR\otimes
id)\DR ~.\label{Dprop2}
\ena
An element $\om$ of
$\Ga$ is said to be {\sl left invariant} if
\eq
\DL (\om) = I \otimes \om \label{linvom}
\en
\noi and {\sl right invariant} if
\eq
\DR (\om) = \om \otimes I \label{rinvom}
\en
\sk
We have seen that any $\rho$ is of the form $\rho=\om^ia_i.$ We have
that the $\om^i$ are left invariant and form a basis of $\invG$, the
linear subspace of all left invariant elements of $\Ga$.
Relation (\ref{braket}) tells us that $\invG$ and $\invX$ are dual
vector spaces.\hfill \\
\sk
\noi iii) There exists an {\sl adjoint representation} $\M{j}{i}$ of the
Hopf algebra, defined by the right action on the $\om^i$:
\eq
\DR (\om^i) = \om^j \otimes \M{j}{i}~;~~~\M{j}{i} \in A~. \label{adjoint}
\en

The co-structures on the $\M{j}{i}$ can be deduced \cite{Wor}:
\begin{eqnarray}
& & \Delta (\M{j}{i}) = \M{j}{l} \otimes \M{l}{i} \label{copM}\\
& & \epsi (\M{j}{i}) = \delta^i_j \label{couM}\\
& & S (\M{i}{j}) =  (M^{-1})_i{}^j
\label{coiM}\end{eqnarray}

The elements $\M{j}{i}$ can be used to
build a right invariant basis of $\Ga$. Indeed the $\eta^i$ defined by
\eq
\eta^i \equiv \om^j S (\M{j}{i}) \label{eta}
\en
are a basis of $\Ga$ (every element of $\Ga$ can be uniquely written
as $\rho = \eta^i b_i$) and their right invariance can be checked
directly .

Moreover, from (\ref{coiM}), using (\ref{eta})
and (\ref{aom})
one can prove the relation
\eq
 (a * \f{j}{i})\M{j}{k}=\M{i}{j}(\f{k}{j} * a)  \label{propM}
\en
\noi with $a* \f{j}{i} \equiv (\f{j}{i} \otimes id) \D(a)$
being the action of the right-invariant operator $\f{j}{i}$
on the function $a$.
\sk

\sect{Construction of the space of Vector Fields.}

In this section we  study the space $\Xi$ of vector fields over Hopf algebras
defining a
left product between elements of $A$ and of $\invX$.
\sk
A generic vector field can be written in the
form $f^i\cdot t_i$ where $\{t_i\}~\,i=1,\,\ldots,n$
is a basis of left invariant
vector fields and $f^i$ are $n$ smooth functions on the group
manifold.
In the commutative case $f^i\cdot t_i=t_i\dia f^i$ i.e. left and right
products (that we have denoted with $\dia$) are  the same, indeed $ (t_i\dia
f^i)(h)\equiv t_i(h)f^i=f^it_i(h) = (f^i\cdot t_i)(h)$.

Let ${t_i} = {\chi_i *}$ be  a basis in  $\invX$ and
let $a^i$, $\,i=1,\,\ldots,n$ be  generic
elements of $A$:
\sk
\noi {\bf Definition}
\eq
\Xi \equiv \{V~|~~ V: A \longrightarrow A ~;~ V = a^i\cdot t_i \}~,
\label{defff}
\en
where the definition of the left product $\cdot$ is given below:
\sk
\noi {\bf Definition}
\eq
\forall a,b \in A  ,\forall t \in {}_{inv}\Xi
{}~~~~~~(a \cdot t) b \equiv a t(b) = a (\chi * b) ~.
\en
The product $\cdot$ has a natural generalization to the whole $\Xi$ :
\eq \begin{array}{rcl}
\cdot~:~& A \times \Xi \longrightarrow & \Xi \nonumber\\
	 & ~ (a,V) ~        \longmapsto     & a\cdot V
\end{array}
{}~~~ \mbox{ where }~~~ \forall b \in A~~~ (a \cdot V)(b) \equiv a V(b) ~.
\en

\sk
\noi It is easy to prove that $( \Xi ,\cdot )$ is a left $A$-module:
\eq
(a + b) \cdot V = a\cdot V + b\cdot V~;~~~
  (ab)\cdot V = a\cdot(b\cdot V)~;~~~
(\lambda a)\cdot V = \lambda a\cdot V
\label{amodule}
\en
with $\lambda\in
\mbox{\boldmath$C$}\,$.

\noi {\sl Note:} {}~ The left product was chosen in the definition
of $\Xi$
so that the symbol $\cdot$ can henceforth be omitted in all formulae. For the
right product $\dia$  --- that we will introduce later ---
we {\em do}\/ have to distinguish for instance
the elements $V\dia(ab) \in \Xi$ and $V(ab) \in A$ because the
vector fields act to the right (by convention).
\sk
$\Xi$ is the analog of the space of derivations on the ring
$C^{\infty}(G)$ of the
smooth functions on the group $G$.
Indeed we have:

\eq V(a+b)=V(a)+ V(b)~~,~~~~~ V(\lambda a) = \lambda
V(a)~~~~\mbox{Linearity}\label{Linearity}
\en
\eq
{}\:~~~V(ab) \equiv (c^i t_i)(ab) = V(a)b + c^i(f_i{}^j*a)t_j(b)
{}~~~~~~\mbox{ Leibniz rule}\label{Leibnizrule}
\en

\noi in the classical case $c^i(f_i{}^j*a)t_j(b) =
aV(b)${}$\,$ (recall $\f{j}{i}=\delta^j_i\epsi\; ;~\:\epsi *b=b).$

\sk
\noindent This last equation can be written as a commutation relation
\eq
V a = V(a) + c^i(f_i{}^j*a)t_j
\en
and in a basis-independent version that
uses the right product $\dia$ (\ref{diadef})\footnote{There is
another basis-independent version of this equation based
on a left $U$-coaction
${}_{U}\Delta(V) \equiv
V^{1'} \otimes V^2$ (that generalizes the coproduct of $\invX$):~
$V a = V^{1'}(a) V^2$.}
\eq
V a = V(a) + V \dia a.\label{Productrule}
\en

We have seen the duality between $\invG$ and $\invX$. We now extend it to
$\Ga$ and $\Xi$, where $\Ga$ is seen as a right $A$-module (not
necessarily a bimodule) and $\Xi$ is our left $A$-module.

\sk

\noi
{\bf Theorem} 1.
There exists a unique map
$$
\langle~~,~~\rangle ~:~~~\Xi \times \Ga \longrightarrow A
$$
\indent such that:
\sk
\noi 1) $ \forall\: V \in \Xi$; the application
\[
 \langle V,~~ \rangle  ~:~~\Ga \longrightarrow A
\]
is a  right $A$-module morphism, i.e. is linear and $\langle
V,\rho a \rangle=\langle V , \rho\rangle a$.

\noi 2) $\forall \rho\in \Ga $; the application
\[ \langle ~~,\rho\rangle ~:~~\Xi \longrightarrow A\]
is a left $A$-module morphism, i.e. is linear and $\langle
b V, \rho\rangle=b\langle V, \rho\rangle $.

\noi 3) Given  $\rho\in\Ga$
\eq
\langle ~~,\rho\rangle=0
  ~\Rightarrow ~ \rho=0~, \label{Duality1}
\en
\sk
\noi where $\langle ~~,\rho\rangle=0$ means
$\langle V,\rho\rangle = 0 ~~\forall V \in \Xi  $.

\noi 4) Given  $V \in \Xi$
\eq
\langle V,~~\rangle=0
{}~\Rightarrow ~ V=0~, \label{Duality2}
\en
\sk
\noi where $\langle V,~~\rangle=0$ means $\langle V,\rho\rangle = 0 ~~\forall
\rho \in \Ga  .$

\noi 5) On $\invX \times \invG$ the bracket $ \langle~~,~~\rangle$ acts as the
one introduced in the previous section.

\sk
\noi {\sl Remark}. Properties 3) and 4) state that
$\Ga$ and $\Xi$ are dual $A$-moduli, in the sense that
they are dual with respect to $A$.

\sk
\noi {\sl Proof}

\noi Properties 1), 2) and 5) uniquely characterize this map . To prove
the existence of such a map we show that the following bracket
\sk
\noi {\bf Definition}
\eq
\langle V,\rho\rangle \equiv V(a_{\al})b_{\al}~,
\en
 where $a_{\al},b_{\al}$ are elements of $A$ such that
$\rho =da_{\al} b_{\al}$,
satisfies 1),2) and 5).
\sk
\noi We first verify that the above definition is well given, that is:
\[
\mbox{Let } \rho =da_{\al} b_{\al}=da'_{\beta} b'_{\beta}~~\mbox{ then }
{}~~V(a_{\al})b_{\al}= V(a'_{\beta})b'_{\beta}
{}~.\]
\noi Indeed,
since \[da_{\al} b_{\al} = da'_{\beta} b'_{\beta}  ~\Leftrightarrow
{}~
       \om^i t_i(a_{\al})b_{\al} = \om^i t_i(a'_{\beta})b'_{\beta}
{}~\Leftrightarrow ~
       t_i(a_{\al})b_{\al} = t_i(a'_{\beta})b'_{\beta}\]\\
{[we used  the uniqueness of the decomposition (\ref{rhoaom})],}
\noi the definition is consistent
 because $$V(a_{\al})b_{\al}=V(a'_{\beta})b'_{\beta}
{}~\Leftrightarrow ~ c^i t_i(a_{\al})b_{\al}=c^i t_i(a'_{\beta})b'_{\beta} $$
where $V=c^i t_i.$
\sk
\noi Property 1) is trivial since $\rho a=(da_{\al} b_{\al})a =
da_{\al}(b_{\al}a).$ \\
\noi Property 2) holds since
\[\langle c V,\rho\rangle = (c V)(a_{\al})b_{\al} =
c(V(a_{\al})b_{\al})=c \langle V,\rho\rangle
{}~.\]

\noi Property 5).
Let $\{\om^i\}$ and $\{t_i\}$ be dual bases in $\invG$ and $\invX$.
Since  $\om^i \in
\Ga ~,~ \om^i = da_{\al} b_{\al}$ for some $a_{\al}$ and $b_{\al}$ in $A$.
We can also write $\om^i =da_{\al} b_{\al}= \om^k t_k(a_{\al})b_{\al}$ ,
so that, due to the uniqueness of the decomposition (\ref{rhoaom}),
we have
\[~~~~~~~~~~~~~~~ t_k (a_{\al})b_{\al}=\delta_k^i I~~~~~~~~~~~~(I
\mbox{ unit of} A);\]

\noi we then obtain  \[\langle t_j,\om^i\rangle=t_j(a_{\al})b_{\al}
=\delta^i{}_j
I~.\]

\noi Property 3).
Let  $\rho= \om^i a_i\in\Ga~.$ \\
If $\langle V , \rho\rangle =0 ~~\forall V\in \Xi$, in particular
$\langle t_j, \rho\rangle = 0 ~~\forall j=1,\ldots ,n$; then
$\mbox{$\langle t_j,\om^i \rangle a_i=0 \Leftrightarrow $} \\ a_j=0$ ,
and therefore $\rho =0~.$
\sk
\noi Property 4).
Let $V = a^i t_i  \in \Xi~.$ \\
If $\langle V ,\rho \rangle =0 ~~\forall\rho\in\Ga$, in particular
$\langle V  , \om^j\rangle = 0 ~~\forall j=1,\ldots ,n$; then
$\mbox{$a^i\langle t_i,\om^j \rangle =0 \Leftrightarrow $} \\ a^j=0$ ,
and therefore $V =0~.$

\cvd
By construction every $V$ is of the form
\[ V= a^i t_i.\]
\indent We can now show the unicity of such a decomposition.
\sk
\noi {\bf Theorem} 2. Any $V \in \Xi$ can be uniquely written in the form
\[V=a^i t_i
\]
\noi {\sl Proof}\\
Let $ V=a^i t_i=a'^i t_i $ then
\[\!\!\!\!\!\!\!\!\!\!\!\!\!\!\!\!\!\!\!
\forall i=1,\,\ldots,n~~~~~~ a^i=  a^j\langle t_j, \om^i\rangle  =
 \langle V ,\om^i \rangle = a'^j \langle t_j ,\om^i \rangle  = a'^i~.
\]
\cvd
\noi Notice that once we know the decomposition of $\rho$ and $V$ in terms of
$\om^i$ and $t_i$, the evaluation of $ \langle~~,~~\rangle$ is trivial:
\[
\langle V,\rho\rangle = \langle b^j t_j , \om^i a_i\rangle = b^j\langle
t_j,\om^i\rangle a_i = b^ia_i~.
\]
Vice versa from the previous theorem $V=\langle V,\om^i\rangle t_i$ and
$\rho = \om_i\langle t_i,\rho\rangle ~.$

\pagebreak

We conclude
this section summarizing the three different ways
of looking at $\Xi$:

\noi (I)${}~{}~{}$   \spz $\Xi$ as the
set of all deformed derivations over $A$ [see (\ref{defff}),
(\ref{Linearity}) --
(\ref{Productrule})].

\noi (II)${}~{}$  \spz $\Xi$ as the left $A$-module freely  generated by
the {\sl elements} $t_i\,, ~~i=1,\,\ldots,n$. The latter is  the
set of all the {\sl formal} products and sums of the type $a^it_i$,
where $a^i$ are generic elements of $A$.
Indeed, by virtue of Theorem 2, the map that associates to each
$V=a^i\cdot t_i $ in $\Xi$ the corresponding element $a^i t_i$
 is an isomorphism between left $A$-moduli.

\noi (III) $\!$\spz  $\Xi$ as $ \Xi '=
\left\{ U ~:~\Ga \longrightarrow A ,
{}~U \mbox{ linear and } U(\rho a)=U(\rho)a ~\forall a\in A\right\}$, i.e.
$\Xi$ as the dual (with respect to $A$) of the space of 1-forms $\Ga$.
The space $\Xi '$ has a trivial left $A$-module structure: $(aU)(\rho)
\equiv aU(\rho)~.$
$\Xi$ and $\Xi '$ are isomorphic left $A$-moduli because of property
(\ref{Duality2}) which states that to each $\langle V,~~\rangle \: :~\Ga
\rightarrow A$
 there
corresponds one and only one
$V$.$[\langle V,~~\rangle=\langle V',~~\rangle\Rightarrow V=V'].$
Every $U\in \Xi '$ is of the form $U = \langle V, ~~\rangle$; more precisely,
if $a^i$ is such that
$ U(\om^i)=a^i$ then $U = \langle a^i t_i,~~~ \rangle~.$

These three ways of looking at $\Xi$ will correspond to different
aspects of the Cartan Calculus: the Lie derivatives $\ell_V$ will
generalize (I), inner derivations $i_V$ will correspond to (III), while
the transformation properties of $\ell_V$ and $i_V$ are governed by (II).

\sk
\sect{Bicovariant Bimodule Structure.}
\sk   \label{secBBS}
In \cite{Wor} the space of $1$-forms is extensively studied.
A right and a left product are introduced between elements of $\Ga$ and of $A$,
and it is known how to obtain a right product from a left one [e.g.
$a\om^i=\om^j(f_j{}^i*a)$, $\om^ia=(Sf_j{}^i*a)\om^j$],
i.e. $\Ga$ is a bimodule over $A$.
\sk
Since the actions $\DL$ and $\DR$ are compatible in the sense that:
\eq
(id \otimes \DR) \DL = (\DL \otimes id) \DR
\en
the bimodule $\Ga$ is called a {\sl bicovariant bimodule}.
In \cite{Wor} it is shown that relations
(\ref{propf1}), (\ref{propf2}), (\ref{aom}), (\ref{propM}),
(\ref{copM}), (\ref{couM})
and (\ref{adjoint}) completely characterize the bicovariant
bimodule $\Ga$.
\sk
In Section 3 we have studied the left product $\cdot$ and we have seen
that $\Xi$ is a left bimodule over $A$ [see (\ref{amodule})].
In this section we introduce a right product and a left and right action
of the Hopf algebra $A$ on $\Xi$. The left and right actions $\DS$ and
$\DD$ are the $q$-analog of the push-forward of tensor fields on a
group manifold. Similarly to $\Ga$ also $\Xi$ is a bicovariant bimodule.
\sk
The construction of the right product on $\Xi$, of the right action $\DD$
and of the left action $\DS$ will be effected  along the lines of Woronowicz'
Theorem 2.5 in \cite{Wor}, whose  statement can be explained  in the
following steps:
\sk
\noi {\bf Theorem} 3.
Consider the {\sl symbols} ${t_i}~ (i=1,\,\ldots,n)$
 and let $\Xi$ be the left $A$-module freely generated by them:

\[
\Xi \equiv \{a^i t_i ~|~~a^i \in A\}
\]
Consider functionals $\FF{i}{j}:~A\longrightarrow \mbox{\boldmath$C$}~$
satisfying [see (\ref{propf1}) and (\ref{propf2})]
\begin{eqnarray}
& & \FF{i}{j} (ab)= \FF{i}{k} (a) \FF{k}{j} (b) \label{propF1}\\
& & \FF{i}{j} (I) = \del^i_j~~~~~~~~~~~~~~~. \label{propF2}
\end{eqnarray}
\indent Introduce a right product via the definition [see (\ref{aom})]
\sk
\noi {\sl Definition}
\eq
(a^i t_i)\dia b \equiv a^i (\FF{i}{j}* b)t_j.  \label{bta}
\en
It is easy to prove that

\noi i) \spz $\Xi$ is  a bimodule over $A$. (A proof of this first
statement as well as of the following ones is contained in \cite{Wor}).

\rightline{$\Box$}
\sk
Introduce a left action (push-forward) of the Hopf algebra $A$ on $\Xi$
\sk
\noi {\sl Definition}
\eq
\DS(a^it_i)\equiv\D(a^i)(I\otimes t_i)\label{DDta}~.
\en
It follows that

\noi ii) \spz $(\Xi, \DS)$ is a left covariant
 bimodule over $A$, that is
\[
 \DS(a V b)=\D (a)\DS(V)\D(b)~;~~~
 (\epsi \otimes id) \DS (V)=V~;~~~
 (\D \otimes id)\DS=(id\otimes\DS)\DS ~.
\]
\rightline{$\Box$}
\sk
Introduce $n^2$ elements of $A$, $\N{i}{j}$,
satisfying [see (\ref{propM}),(\ref{copM}) and
(\ref{couM})]
\begin{eqnarray}
& \N{i}{j} (a * \FF{i}{k})=(\FF{j}{i} * a) \N{k}{i} & \label{propN} \\
 & \Delta (\N{j}{i}) = \N{j}{l} \otimes \N{l}{i} &  \label{copN}\\
 & \epsi (\N{j}{i}) = \delta^i_j~~~~~~, &  \label{couN}
\end{eqnarray}
and introduce a right action $\DD $ such that [see (\ref{adjoint})]
\sk
\noi {\sl Definition}
\eq
\DD(a^i t_i)\equiv \D(a^i) (t_j\otimes \N{j}{i})\label{DS}~.
\en
Then it can be proven that

\noi iii) \spz The elements  [see (\ref{eta})]
\eq
h_i \equiv S^{-1}(\N{j}{i})t_j \label{heich}
\en
are right invariant: $\DD(h_i) = h_i\otimes I$.
Moreover any $V \in \Xi$ can be expressed in a unique way respectively as
$V=a^ih_i$ and as
$V=h_ib^i$, where $a^i, b^i \in A$.

\rightline{${}~~~~~~~~~~~~~~~~~~~~~~~~~\Box$}

\noi iv) \spz $ (\Xi, \DD)$ is a right covariant bimodule over $A$,
that is
\[
 \DD(a V b)=\D (a)\DD(V)\D(b)~;~~~ (id \otimes \epsi) \DD (V)=V~;~~~
 (id \otimes \D)\DD =(\DD\otimes id)\DD~.
\]{}
\rightline{$\Box$}

\noi v) \spz The left and right covariant bimodule $(\Xi,\DS,\DD) $ is a
bicovariant
bimodule, that is left and
right actions are compatible:
\[
(id\otimes\DD)\DS=(\DS\otimes id)\DD~.
\]
\cvd

In the previous section we have seen [remark (II)] that the space of vector
fields
$\Xi$ is the free left $A$-module generated by the symbols $t_i$, so
that the above theorem applies to our case.

There are many bimodule structures (i.e. choices of $\FF{i}{j}$) $\Xi$ can
be endowed with.
Using the fact that $\Xi$ is dual to $\Ga$ we request compatibility with
the $\Ga$ bimodule.\\
In the commutative case $\langle ft_j,\om^i\rangle=
\langle t_j\dia f,\om^i\rangle=
\langle t_j,f\om^i\rangle=\langle t_j,\om^if\rangle.$\\
In the quantum case we know that $\langle t_j, \om^ia\rangle =
\langle a t_j,\om^i\rangle$ and we require
\eq
\langle t_j\dia a,\om^i\rangle =
 \langle t_j,a\om^i\rangle~; \label {wat}\en
this condition uniquely determines the bimodule structure of $\Xi$.
Indeed we have
\eq
t_i\dia a = (\f{j}{i} * a~) t_j \label{diadef}
\en
i.e.
\eq
\FF{i}{j} \equiv \f{j}{i}.
\en
This is in fact the bimodule structure needed for the
product rule (\ref{Productrule}).
So far $\Xi$ has  a bimodule structure.
$\Xi$ becomes a left covariant bimodule if we define $\DS$ as in
(\ref{DDta}) so that $t_i$ are left invariant vector fields. The
following consideration fixes the right covariance of $\Xi$:
For reasons of symmetry let us try to impose
\eq
\langle t_i , \om^j \rangle = \delta_i^j = \langle h_i , \eta^j \rangle,
\en
where $t_i$ and $\om^j$ are left-invariant and $h_i$ and $\eta^j$ are
the canonically associated right-invariant objects; see
(\ref{heich}) and (\ref{eta}).
As can easily be seen this condition is equivalent to the following
definition of the $n^2$ elements $\N{l}{k} \in A$
\sk
\noi {\bf Definition}
\eq
\N{l}{k}=S(\M{k}{l})~.   \label{NeqSM}
\en
\sk
\noi {\bf Theorem} 4.  The $\N{l}{k}$ elements defined above satisfy relations
(\ref{couN}),
(\ref{copN}) and
(\ref{propN}):
$$
1)~ \epsi (\N{j}{i}) = \delta^i_j
{}~~~2)~ \Delta (\N{j}{i}) = \N{j}{l} \otimes \N{l}{i}
{}~~~3)~ \N{i}{k} (a * \FF{i}{j})=(\FF{k}{i} * a) \N{j}{i}
$$
\sk
\noi {\sl Proof}
\sk
\noi $1)~$ is trivial.
\sk
\noi $2)~$ use $N^i{}_j = SM_j{}^i$ and $
\D\circ S = \sigma \circ (S \otimes S)\circ \D$,
where $\sigma_A$ is the flip map in $A\otimes A$.
\sk
\noi $3)~$ Using $N^i{}_j = SM_j{}^i$ and $O_i{}^j = f_i{}^j$
in (\ref{propM}) gives
$$
(a * O_k{}^i) S^{-1}N^j{}_i = S^{-1}N^i{}_k (O_i{}^j * a).
$$
Multiplying with ``$N$'' on both sides we obtain relation 3).
 \sk \cvd
In equations (\ref{donf}) and (\ref{dachi})
elements of $\Xi$ and $\Gamma$ make a
joint appearance. To be still able to talk about transformation
properties of such expressions we need to combine the previously
introduced actions into one object, $\Delta_A$, simply by putting
$\Delta_A \equiv {}_\Xi\Delta$ on $\Xi$ and $\Delta_A \equiv
{}_\Gamma\Delta$ on $\Gamma$ and requiring $\Delta_A$ to be an
algebra homomorphism. From this definition we
get the following important corollary:
\sk
\noindent {\bf Corollary.} The expression $\omega^i t_i$ in (\ref{donf})
is invariant in the sense that
$$
\Delta_A(\omega^i t_i) = {}_\Gamma\Delta(\omega^i){}_\Xi\Delta(t_i) =
\omega^j t_k \otimes M_j{}^i N^k{}_i = \omega^i t_i \otimes 1.
$$
Similar statements apply to ${}_A\Delta$.
\sk
The left-invariant vector fields were characterized in (\ref{leftXi})
through their left action $t_i(a) = \chi_i * a$ on functions. It would
be nice if a similar equation were true for the right-invariant vector
fields $h_i$; indeed we have
\sk
\noi {\bf Theorem} 5.\spz  $ h_i(a) \equiv
S^{-1}(N^j{}_i) t_j(a) = a * \chi_i~~~~~\forall a \in A.  $
\sk
\noi {\sl Proof}\\
\noi From (\ref{donf}) we have
$$
\om_a \equiv P da \equiv d(a_2) S^{-1}a_1 = \om^i\chi_i(a)
$$
and in particular $\om_{x^i} = \om^i$. Using this equation
we can rewrite the definition of $M_i{}^k$, (\ref{adjoint}),
as $M_i{}^k = \chi_i(x^k_2) x^k_3 S^{-1}x^k_1$, or
$$
M_i{}^j \chi_j(x^k) = \chi_i(x^k_2) x^k_3 S^{-1}x^k_1.
$$
We now show that this equation is still true for
an arbitrary function $a \in A$ in place of $x^k$. Indeed from (\ref{axri})
any function can be decomposed into a part that is spanned by the $x^i$
and a part that is contained in
$R \oplus \{I\}$ \ie is annihilated
by the $\chi_j$. But $R$ and $\{I\}$ are Ad-invariant so the part of
$a$ that is contained in $R \oplus \{I\}$ contributes to neither side of the
equation in question and we have in fact:
$$
M_i{}^j \chi_j(a) = \chi_i(a_2) a_3 S^{-1}a_1.
$$
Definition (\ref{NeqSM}) and a slight rearrangement
(take $a=b_2$ and multiply by $b_1$) finally gives
$$
S^{-1}N^j{}_i \: a_1 \chi_j(a_2) = \chi_i(a_1) a_2
$$
and proves the theorem.
\sk \cvd
\noi {\sl Remark.}
\noi Statement 3) of Theorem~4 is in fact a consequence of
Theorem 5. Both are a consequence of a more general Hopf algebra
relation ({\em e.g.}\ eqn. 4.67 of \cite{S1}) that is valid
in the semi-direct product algebra of $A$ and $A^*$.
\sk
Now that we have all the ingredients, the construction of the bicovariant
bimodule
$\Xi$ is easy and straightforward.
For example $\DD$ is given in formula (\ref{DS}).
\sk
\noi We can then conclude that $(\Xi ,\DS ,\/\DD)$ is a bicovariant
bimodule.
\sk
Notice that, since  Theorem 3  completely characterizes a bicovariant
bimodule all the formulas containing the symbols $\f{j}{i}$ or
$\M{k}{l}$ or elements of $\Ga$ are still valid under the
substitutions $\f{j}{i}\rightarrow S^{-1}\FF{i}{j},~~\M{k}{l}
\rightarrow \N{l}{k} $ and $ \Ga\rightarrow \Xi$.
\sk

\sect{Tensor fields}

The construction completed for vector fields is readily generalized to
$p$-times
contravariant tensor fields.

We define $\Xi\otimes\Xi$ to be the space of all elements that
can be written as finite sums of the kind $\sum_i V_i\otimes V'_i $
with $V_i,V'_i\in\Xi$. The tensor product (in the algebra $A{)}$
between $V_i$ and $V'_i$ has the following properties:
\begin{equation}
 V\dia a\otimes
V'=V\otimes aV' ,~ a(V\otimes V') =(aV)\otimes V' \mbox{ and } (V\otimes
V')\dia a=V\otimes(V'\dia a) \label{tensorA}
\end{equation}
 so that $\Xi\otimes\Xi $ is naturally a bimodule
over $A$.\\
Left and right actions on $\Xi\otimes\Xi$ are defined by:
\eq
{}_A\Delta (V \otimes V')\equiv   V_{1}   {V'} _{1} \otimes   V_l \otimes
  {V'}_l,~~~{}_A\Delta: \Xi \otimes \Xi \rightarrow A\otimes\Xi\otimes\Xi
\label{DSXiXi}
\en
\eq
\Delta_A (  V \otimes   {V'})\equiv   V_r \otimes   {V'}_r \otimes   V_1
  {V'}_1,~~~\Delta_A: \Xi \otimes \Xi \rightarrow \Xi\otimes\Xi\otimes A
\label{DDXiXi}
\en
\noi where $  V_{1} \otimes V_l$ and $V_r \otimes V_1$ are defined by
\eq
\DS (  V) =   V_{1} \otimes   V_l~~~  \in A \otimes \Xi~,
\en
\eq
\DD (  V) =   V_r \otimes   V_1~~~  \in \Xi \otimes A~.
\en
\noi More generally, we can introduce the actions
on $\Xi^{\otimes p}\equiv\underbrace{\Xi \otimes \Xi \otimes \cdots
\otimes \Xi}_{\mbox{$p$-times}}$ as
\[
{}_A\Delta (  V \otimes   {V'} \otimes \cdots \otimes   {V''})\equiv
  V_{1}   {V'}_{1} \cdots   {V''}_{1} \otimes   V_l \otimes
  {V'}_l\otimes \cdots \otimes   {V''}_l
\]
\eq
{}_A\Delta~:~~ \Xi^{\otimes p} \longrightarrow
A\otimes\Xi^{\otimes p}~;
\label{DSXiXiXi}
\en
\[
\Delta_A (  V \otimes   {V'} \otimes \cdots \otimes   {V''})\equiv
  V_r \otimes   {V'}_r \otimes \cdots \otimes   {V''}_r \otimes   V_1
  {V'}_1 \cdots   {V''}_1
\]
\eq
\Delta_A~:~~ \Xi^{\otimes p} \longrightarrow
\Xi^{\otimes p}\otimes A~~.
\label{DDXiXiXi}
\en

\noi Left invariance on $\Xi\otimes\Xi$ is naturally defined as
${}_A\Delta (  V \otimes   {V'}) = I \otimes   V \otimes   {V'}$ (similar
definition for right invariance), so that for example $t_i \otimes
t_j$ is left invariant, and is in fact a left invariant basis for $\Xi
\otimes \Xi$: each element can be written as $a^{ij} t_i\otimes t_j$ in a
unique way.

 It is not difficult to show that $\Xi \otimes \Xi$ is a bicovariant bimodule.
In the same way also
$(\Xi^{\otimes p},{}_A\Delta,\Delta_A)$ is a bicovariant bimodule.
An analog procedure, using $\DL$ and $\DR$ instead of $\DS$ and $\DD$,
applies also to
$\Ga^{\otimes n}$ the $n$-times tensor product of $1$-forms.
\sk
Any element $v\in \Xi^{\otimes p}$ can be written as
$v=b^{i_1\ldots i_p} t_{i_1}\otimes\ldots t_{i_p}$ in a unique way,
similarly any element $\tau\in\Ga^{\otimes n}$ can be written as
$\tau=\om^{i_1}\otimes\ldots \om^{i_n} a_{i_1\ldots i_n}$ in a unique way.
It is now possible to generalize the previous bracket
$\langle~~,~~\rangle\,
:\;\Xi\times\Ga\rightarrow A\:$ to $\Xi^{\otimes p}$ and $\Ga^{\otimes n}$.
For $p > n$: $\langle v , \tau \rangle = 0$; for $p \leq n$:
\begin{equation}
\begin{array}{rccll}
\langle~~,~~\rangle ~:~~&\Xi^{\otimes p} \times \Ga^{\otimes n}&
\longrightarrow & \Ga^{\otimes n - p} &{} \\
&(v,\tau) &\longmapsto
&\langle v,\tau\rangle &=b^{j_1\ldots j_p}\langle
t_{j_1}\otimes \ldots t_{j_p}\, ,\,\om^{i_1}\otimes
\ldots \om^{i_n}\rangle a_{i_1\ldots i_n}\\
&{}&{}&{}&=b^{i_p\ldots i_1}\om^{i_{p+1}}\otimes \ldots
\om^{i_n} a_{{i_1}\ldots i_n}
\end{array}\label{genbrack}
\end{equation}
where $\Ga^{\otimes 0}
\equiv A$, $\Ga^{\otimes 1} \equiv \Ga$ and we have defined
\eq
\begin{array}{rcl}
\langle t_{j_p}\otimes \ldots t_{j_1}\, ,\,\om^{i_1}\otimes
\ldots \om^{i_n}\rangle
& \equiv & \langle t_{j_p} \otimes
\ldots t_{j_1}, \om^{i_1} \otimes \ldots
\om^{i_p} \rangle \om^{i_{p+1}} \otimes \ldots \om^{i_n} \\
& \equiv & \langle t_{j_1},\om^{i_1}\rangle
\ldots \langle t_{j_p},\om^{i_p}\rangle
\om^{i_{p+1}}\otimes \ldots \om^{i_n}\\
& = & \delta^{i_1}_{j_1} \ldots \delta^{i_p}_{j_p}
\om^{i_{p+1}}\otimes \ldots \om^{i_n}
\end{array}
\label{mirrorcoup}
\en
(with the obvious convention for the special cases $p = n-1$ and $p = n$).
Using definition (\ref{mirrorcoup}) it is easy to prove that
\eq
\langle v \dia a,\tau\rangle=\langle v,a \tau\rangle \label{pass}~,
\en
namely
$$
\begin{array}{rcl}
\langle t_{j_1}\otimes\ldots t_{j_p} \dia a ,
\om^{i_1}\otimes\ldots\om^{i_n}\rangle
& = & (\f{k_1}{j_1}*\ldots \f{k_p}{j_p}*a)\langle
t_{k_1}\otimes\ldots t_{k_p} ,
\om^{i_1}\otimes\ldots\om^{i_n}\rangle\\
& = & (\f{i_p}{j_1}*\ldots \f{i_1}{j_p}*a)
\om^{i_{p+1}}\otimes \ldots \om^{i_n}
\end{array}
$$
$$
\begin{array}{rcl}
\langle t_{j_1}\otimes\ldots t_{j_p},a \om^{i_1}\otimes\ldots\om^{i_n}\rangle
& = & \langle t_{j_1}\otimes\ldots t_{j_p},
\om^{l_1}\otimes\ldots\om^{l_n}\rangle
(\f{i_n}{l_n}*\ldots \f{i_1}{l_1}*a)\\
& = & \delta_{j_p}^{l_1}\ldots\delta_{j_1}^{l_p}
\om^{l_{p+1}}\otimes \ldots \om^{l_n} (\f{i_n}{l_n}*\ldots \f{i_1}{l_1}*a)\\
& = & (\f{i_p}{j_1}*\ldots \f{i_1}{j_p}*a)
\om^{i_{p+1}}\otimes \ldots \om^{i_n}
\end{array}
$$
and (\ref{pass}) is verified if and {\sl only if}
(\ref{mirrorcoup}) holds.

Therefore we have also shown that definition (\ref{mirrorcoup}) is the
only one compatible with property (\ref{pass}), i.e. property
(\ref{pass}) uniquely determines the coupling between $\Xi^{\otimes}$
and $\Ga^{\otimes}.$

It is also easy to prove that
the bracket $\langle~~,~~\rangle$ extends
to $\Xi^{\otimes p}$  and $\Ga^{\otimes p}$ the duality between
$\Xi$ and $\Ga$.
\sk
More generally we can define $\Xi^{\otimes}\equiv
A\oplus\Xi\oplus\Xi^{\otimes 2}\oplus\Xi^{\otimes 3}\ldots $
 to be the algebra
of contravariant tensor fields (and
$\Ga^{\otimes}$ that of covariant tensor fields).

The actions ${}_A\Delta $ and $\Delta_A $ have a natural generalization to
$\Xi^{\otimes}$ so that we can conclude that
$(\Xi^{\otimes},{}_A\Delta ,\Delta_A )$
is a bicovariant graded algebra, the graded algebra of tensor fields
over the ring ``of functions on the group'' $A$, with the left and
right ``push-forward'' ${}_A\Delta $ and $\Delta_A $.

Before we can introduce a Lie derivative and an inner derivation we need
a wedge product of forms, {\em i.e.} we have to
briefly discuss antisymmetrized
covariant tensor fields.

\subsection{Bicovariant graded algebras}

A general method to construct new tensor fields starting from
$\Gamma^\otimes$ is to quotient the bicovariant graded algebra
$\Gamma^\otimes$ with a graded ideal $S$. In order for $\Gamma^\otimes /
S$ to be bicovariant, the ideal $S$ has also to be biinvariant (i.e.
a left and right sub-comodule). For example the ideal can be generated
by the elements $S^{kl}_{ij} \om^i \otimes \om^j$
\footnote{Or, what is the same, one imposes
{\em relations}\ $S^{kl}_{ij} \om^i \wedge \om^j = 0$.}
and a sufficient
condition for the biinvariance is
$S_{ij}^{kl} M_k{}^m M_l{}^n = M_i{}^k M_j{}^l S_{kl}^{mn}$,
where ${}_\Xi\Delta \om^i = \om^j \otimes M_j{}^i$.

Particular cases of such a quotient procedure are when $S$ is the
kernel of a graded bimodule homomorphism $A:\Gamma^\otimes \rightarrow
\Gamma^\otimes$ then $\Gamma^\otimes / {\rm ker}(A)$ is isomorphic
to the image of $A$ and therefore is a biinvariant graded subalgebra
of $\Gamma^\otimes$.
This is what we require for a wedge product construction. The wedge
product of two forms is given by
\begin{equation}
\begin{array}{c}
\om^i \wedge \om^j = A(\om^i \otimes \om^j)
= \om^k \otimes \om^l W_{kl}{}^{ij} \in \Ga^{\wedge 2} \subset
\Ga^{\otimes 2} \\
a\om^i\wedge\om^j\equiv a(\om^i\wedge\om^j)~;~~\om^i\wedge\om^jb\equiv
(\om^i\wedge\om^j)b
\label{Www}
\end{array}
\end{equation}
Since the homomorphism $A$ is required to commute with the coactions
${}_A\Delta$ and $\Delta_A$
the new tensor transforms covariantly according
to its index structure:
\begin{equation} \Delta_A (\om^i \wedge \om^j)
\equiv \Delta_A (\om^k \otimes \om^l ) W_{kl}{}^{ij}
= \om^k \wedge \om^l \otimes M_k{}^i M_l{}^j
\end{equation}
thus imposing
\begin{equation}
W_{ij}{}^{kl} M_k{}^m M_l{}^n = M_i{}^k M_j{}^l W_{kl}{}^{mn} \label{WMM}
\end{equation}
on $W$ as a necessary condition. In short tensor product notation:
$\om_1 \wedge \om_2 = \om_1 \otimes \om_2 W_{12} \Rightarrow
W_{12} M_1 M_2 = M_1 M_2 W_{12}$.
Let $\hat{\sigma}_{ij}{}^{kl} =
f_i{}^l(M_j{}^k)$. Due to (\ref{propM}) any function
$W(1, \hat{\sigma})$ will
satisfy (\ref{WMM}). We will in fact take $W = 1 - \hat{\sigma}$.
For this choice and for examples
of $n > 2$ relations we would like to refer the reader to
\cite{Wor,AC,SWZ3,S1}; here
we will just quote the results:
generalizing  equation (\ref{Www}),
wedge products of $n$ forms are again expressed in terms
of tensor fields:
\begin{equation}
\om_1 \wedge \ldots \wedge \om_n = \om_1 \otimes \ldots \otimes \om_n
W_{1\ldots n}.\label{wedge1}
\end{equation}
The numerical coefficients $W_{1\ldots n}$ are given through a recursion
relation
\begin{equation}
W_{1\ldots n} = W_{2\ldots n} {\cal I}_{1\ldots n}, \label{wedge2}
\end{equation}
where
\begin{equation}
{\cal I}_{1\ldots n} = 1 - \hat{\sigma}_{12} +
\hat{\sigma}_{12} \hat{\sigma}_{23} - \ldots
-(-1)^n \hat{\sigma}_{12} \hat{\sigma}_{23} \cdots \hat{\sigma}_{n-1,n}
\label{wedge3}
\end{equation}
and $W_i{}^j = {\cal I}_i{}^j = \delta_i^j$.
As is easily seen, ${\cal I}$ has the following
decomposition property that we will use
in the next chapter:
\begin{equation}
{\cal I}_{1\ldots n} = {\cal I}_{1\ldots k}
+(-1)^k \hat{\sigma}_{12} \hat{\sigma}_{23} \cdots \hat{\sigma}_{k,k+1}
{\cal I}_{k+1\ldots n}. \label{wedge4}
\end{equation}
The first order differential operator $d$ can be extended in a unique way
to an exterior differential $d$ mapping n-forms in (n+1)-forms,
satisfying the undeformed
Leibniz rule and such that $d^2=0$ \cite{Wor}.
\sk
\noi It is also straightforward to define antisymmetrized contravariant
vector fields in this way.

\sect{Contraction operator and Lie derivative}
In this last section we propose a definition of Lie derivative along a
generic vector field.
We start with the introduction of  the contraction operator $i_V $, also
called inner derivation, with
$V\in\Xi$ based on equation (\ref{genbrack}).
\sk
\noi For a generic vector field  $V=b^j t_j$ we define:
\nopagebreak
\sk
\nopagebreak
\noi {\bf Definition}
$$
i_V(\vartheta) = \langle V , \vartheta \rangle
\qquad \vartheta ~ \mbox{generic form.}
$$
\sk
\noi {\bf Theorem} 6. The contraction operator satisfies the
following properties:\\
$\vartheta, \vartheta'$ generic forms; $a, a_{i_1\ldots i_n} \in A$;
$V=b^i t_i$; $\lambda \in$ {\bf C}
\begin{enumerate}
\item[a)] $i_V(\vartheta) = i_{b^j t_j}(\vartheta) =
b^j i_{t_j}(\vartheta)$
\item[b)] $i_{V}(a)=0$
\item[c)] $i_{V} ( \omega^{j} )= b^j$
\item[d)] $i_{V} (\omega^{i_{1} }\wedge \ldots  \omega^{i_{n}}
      a_{ i_{1} \ldots i_{n} }  )
      =
      \begin{array}[t]{l}
	i_{V} (\omega^{i_{1}} \wedge \ldots \omega^{i_{s}})
	\wedge (\omega^{i_{s+1} } \wedge  \ldots \omega^{i_{n}}
	a_{ i_{1} \ldots i_{n} }) + \\
	(-1)^{s} b^i
	\f{j}{i} * (\omega^{i_{1} } \wedge \ldots\omega^{i_{s}} )
	\wedge i_{t_{j}} (\omega^{i_{s+1} } \wedge  \ldots\omega^{i_{n}}
	a_{ i_{1} \ldots i_{n} })
      \end{array}$
\item[e)] $i_{V} (a \vartheta) = b^i (f_i{}^j * a)\, i_{t_j}(\vartheta)$
\item[f)] $i_{V} (\vartheta a +\vartheta') =
      i_{V} ( \vartheta ) a + i_{V}(\vartheta')$
\item[g)] $i_{\lambda V} = \lambda i_{V}$
\end{enumerate}
\sk
\noi {\sl Proof}.\\
Properties a), b), c), f), g) are direct consequences of (\ref{genbrack}),
e) follows from (\ref{pass}). To proof d) we have to use the definition
of the wedge product (\ref{wedge1} -- \ref{wedge3}): first we note that
(in tensor product notation)
$$\begin{array}{rcl}
i_V(\om_1 \wedge \ldots \wedge \om_n)
& = & i_V(\om_1 \otimes \ldots \otimes \om_n W_{1\ldots n})\\
& = & i_V(\om_1)\om_2 \otimes \ldots \otimes \om_n W_{1 \ldots n}\\
& = & i_V(\om_1)\om_2 \wedge \ldots \wedge \om_n {\cal I}_{1 \ldots n},
\end{array}$$
where we have used (\ref{wedge2}) in the last step --- in index notation:
$$i_{t_i}(\om^{i_1}\wedge \ldots \wedge \om^{i_n})
= \om^{j_2}\wedge \ldots \wedge \om^{j_n} {\cal I}_{i j_2 \ldots j_n}^{i_1
\ldots i_n};$$
next we can show
$$\begin{array}{rcl}
\f{j}{i} * (\omega^{i_{1} } \wedge \ldots\wedge\omega^{i_{s}} )
& = & \omega^{k_{1} } \wedge \ldots\wedge\omega^{k_{s}}
      \f{j}{i}(M_{k_1}{}^{i_1} \cdots M_{k_s}{}^{i_s})\\
& = & \omega^{k_{2} } \wedge \ldots\wedge\omega^{k_{s+1}}
      \hat\sigma_{i k_2}{}^{i_1 l_2} \hat\sigma_{l_2 k_3}{}^{i_2 l_3} \ldots
      \hat\sigma_{l_s k_{s+1}}{}^{i_s j};
\end{array}$$
finally we utilize the decomposition
property (\ref{wedge4}) and associativity of the wedge product
(in tensor product notation)
$$\begin{array}{rcl}
i_{V} (\om_1 \wedge &&\!\!\! \ldots \wedge \om_n) ~=\\
 = && i_V(\om_1)\om_2 \wedge  \ldots \om_n {\cal I}_{1 \ldots n}\\
 = && i_V(\om_1)(\om_2 \wedge \ldots \om_s)
       \wedge (\om_{s+1}\wedge\ldots \om_n)({\cal I}_{1 \ldots s} +
       (-1)^s \hat\sigma_{12}\cdots\hat\sigma_{s,s+1}{\cal I}_{s+1 \ldots n})\\
 = && \begin{array}[t]{l}
      i_V(\om_1\wedge \ldots \wedge \om_s)\wedge(\om_{s+1}\wedge\ldots
      \wedge\om_n)\\
      +b^i i_{t_i}(\om_1)(\om_2 \wedge \ldots \om_{s+1})\wedge
      (\om_{s+2}\wedge\ldots \om_n
      (-1)^s \hat\sigma_{12}\cdots\hat\sigma_{s,s+1}{\cal I}_{s+1 \ldots n})
      \end{array}\\
 = && \begin{array}[t]{l}
      i_V(\om_1\wedge \ldots \wedge \om_s)\wedge\om_{s+1}\wedge\ldots
      \wedge\om_n\\
      +(-1)^s b^i f_i{}^j * (\om_1 \wedge \ldots \wedge \om_s)
      \wedge i_{t_j}(\om_{s+1}\wedge\ldots\wedge\om_n).
      \end{array}
\end{array}$$
With property f) this proves d).\\
\cvd
\sk
\noi {\sl Remark 1}. A slight generalization of property d) for two
generic forms $\vartheta$ and $\vartheta'$ is also true [use e)]:
\eq
i_V(\vartheta \wedge \vartheta') =
i_V\vartheta \wedge \vartheta'
+ (-1)^{{\rm deg}
(\vartheta)} b^i f_i{}^j * \vartheta \wedge
i_{t_j} \vartheta'.\label{itwotheta}
\en
\noi {\sl Remark 2}. Property e) reduces in
the commutative case to the familiar formula:
\[
i_V(h\vart) = hi_V\vart ~.
\]
\sk
\noi The Lie derivative along left invariant vector fields
is given by \cite{AC}:
\eq\ell_t(\tau) \equiv (id \otimes \chi)\DR (\tau) = \chi * \tau \;\;
\;\;\;\;\; \ell_t\; : \;\; \Gamma^{\otimes n} \longrightarrow
 \Gamma^{\otimes n}
\;\;\;
\en
where $\chi \in \invX$ is such that $\chi * =t.$
On $\Gamma^\wedge$ it can be proven that \cite{SWZ1,AC,S2} :
\eq
\ell_{t_i}= i_{t_i}  d +d  i_{t_i}~.\label{dit}
\en
It is then natural to introduce the Lie derivative along a generic
vector field $V$ through the following
\sk
\noi {\bf Definition}
\eq
\ell_{V}= i_{V}  d +d  i_{V}   \label{ellV}
\en
\sk
\noi {\bf Theorem} 7. The Lie derivative satisfies the following properties:
\begin{enumerate}
\item[1)] $\ell_Va=V(a)$
\item[2)] $\ell_Vd\vart=d\ell_V\vart$
\item[3)] $\ell_V(\lambda\vart+\vart')=\lambda\ell_V \vart +\ell_V \vart'$
\item[4)] $\ell_{b V}\vart = b \ell_V\vart + db \wedge i_V\vart,$
where $\vart$ is a generic form
\item [5)] $\ell_V(\vart \wedge \vart') = \ell_V\vart\wedge\vart'
	   + b^i (f_i{}^j * \vart) \wedge \ell_{t_j}\vart'
	   +(-1)^{{\rm deg}(\vartheta)}db^i\wedge (f_i{}^j*\vart)\wedge
	   i_{t_j}(\vart')$,

where $\vart$, $\vart'$ are generic forms and $V = b^i i_{t_i}$.
\end{enumerate}
\noi {\sl Proof.}\\
\noi Properties 1), 2), 3) and 4) follow directly from the definition
(\ref{ellV}).\\
Property 5) is also a consequence of definition (\ref{ellV}); the
proof makes use of (\ref{itwotheta}) and
$d(\f{k}{j}*\vart)=\f{k}{j}*d\vart$.\\
\cvd

\section{Commutation relations}

The careful choice of conventions (see introduction) allows that all
relations can be summarized as deformed commutation relations
\cite{SWZ3} of
right-acting operators
($\vart$: $p$-form, $~~ V = b^i t_i$)
\begin{eqnarray}
d \vart & = & d(\vart) + (-1)^p \vart d\\
i_V \vart & = & i_V(\vart) + (-1)^p b^i (f_i{}^j * \vart) i_{t_j}\\
\ell_V \vart & = & \ell_V(\vart) + b^i (f_i{}^j * \vart) \ell_{t_j}
+(-1)^p db^i(f_i{}^j*\vart) i_{t_j} \label{7.7}
\end{eqnarray}
and their actions ($a \in A$, rest as above) as given by:
\eq
i_V(a)  =  0~;~~
i_V(\om^j) = b^j~;~~
\ell_V(a)  =  V(a) = b^i \chi_i * a ~;~~ \ell_V\om^j=b^i\chi_i *\om^j+db^j~.
\en
These operators form a  graded quantum Lie algebra
\begin{eqnarray}
\{ d , d \} & = & 0\\
\left[ d , \ell_V \right] & = & 0\\
\{ d , i_V \} & = & \ell_V
\end{eqnarray}
which is supplemented by two more
relations\footnote{Relations for the inner derivations
can be constructed along the lines of section~5.1 .}
\begin{eqnarray}
& &\left[ \ell_{t_i} ,
\ell_{t_k} \right]_{\hat B} =
\ell_{[t_i,t_k]_{\hat B}}
		    =f_i{}^l{}_k  \ell_{t_l} \\
& &\left[ \ell_{t_i} , i_{t_k} \right]_{\hat B} = i_{[t_i,t_k]_{\hat B}}
		    =  f_i{}^l{}_k  i_{t_l}
\end{eqnarray}
where $\hat B_{ik}^{rs} \equiv
(\widehat{\sigma^{-1}})_{ik}^{rs} = f_i{}^s(N^r{}_k)$,
$\left[ \ell_{t_i} ,
\ell_{t_k} \right]_{\hat B} \equiv
\ell_{t_i}\ell_{t_k}-\hat B_{ik}^{rs}\ell_{t_r}\ell_{t_s}$
and $f_i{}^l{}_k = \chi_i(N^l{}_k).$\footnote{In these notations
the Cartan-Maurer formula reads:\\ $d\om^i\otimes
\chi_i =\om^i\otimes\om^j\otimes \chi_j\chi_i=
-\om^a\wedge\om^b{\hat\sigma}^{pq}_{ab}[\chi_q,\chi_p]_{\hat B}$.
In particular $d\om^i=\om^a\otimes\om^b\chi_a(M_b{}^i)$.}
\sk

Using the right product (\ref{diadef}) the commutation relations
can be restated in basis-independent form as
\eq
V \vart = V(\vart) + V \dia \vart. \label{basind}
\en
Here and in the following $V$ stands for any of the first-order
operators $V=b^i t_i$ (for $\vart \in A$ only) or $i_V$, $\ell_{t_i}$ and $d$.
For example:
$$
t_i a = t_i(a) + t_i \dia a = \chi * a + (f_i{}^j * a) t_j ~;~~
d\vart=d(\vart) + d\dia \vart=d(\vart)+(-1)^{{\rm deg}(\vartheta)}\vart d~.
$$
For $\vart$ a $p$-form one must add the following rules to the definition
(\ref{diadef}) of the right product to take care of the grading:
\eq
d \dia \vart
= (-1)^p \vart d ~;~~ i_V \dia \vart = (-1)^p b^i (f_i{}^j * \vart) i_{t_j}~.
\en
Note that $V$ appears twice in the r.h.s. of (\ref{basind}). These are the two
parts of the (classical) coproduct
\eq
\delta(V) = V \otimes 1 + 1 \otimes V.
\en
In a hopefully self-explanatory notation we rewrite (\ref{basind})
as
\eq
V \vart = \delta(V) \dia \vart
\en
so that the role of $\delta(V)$ is emphasized.

We now extend the definition of $\delta$ to composition of
operators, say $V \cdot V' \cdots V''$:
\eq
\forall \vart\in \Gamma^{\wedge}~~~~~~~\delta(V \cdot V' \cdots V'') \dia
\vart\equiv V \cdot V' \cdots V'' \vart \label{dvt} ~.
\en
{}From equation (\ref{propM}) one can prove that
\eq
f_i{}^k*[i_{t_j}(\vart)]=
\hat B^{rs}_{ij}i_{t_r}[f_s{}^k*\vart] \label{fibif}
\en
the  following theorem is then easily derived:
\sk
\noi {\bf Theorem 8.}\\
The operators $i_V$, $\ell_{t_i}$, $d$ and $V$ ($V=b^i t_i$ is only
defined on $A$) form a braided tensor
algebra.
The braidings are given by:
$$(1 \otimes a^ii_{t_i})(b^ji_{t_j} \otimes 1) =
 - a^i(f_i{}^h*b^j)\hat B^{rs}_{hj}(i_{t_r} \otimes i_{t_s})~;~~
(1\otimes d)(d\otimes 1)=-d\otimes d $$
\eq
(1\otimes d)(i_V\otimes 1)=-i_V\otimes d~;~~
(1\otimes i_V)(d\otimes 1)=-d\otimes i_V ~.\label{braid}
\en
Using (\ref{dit}) we also get all the other braidings, e.g.
\eq
(1 \otimes a^i\ell_{t_i})(b^j\ell_{t_j} \otimes 1) =
a^i(f_i{}^h*b^j)\hat B^{rs}_{hj}(\ell_{t_r} \otimes \ell_{t_s}).
\en
{}From definition (\ref{dvt}) and the above theorem
it follows:
\sk
\noi {\bf Theorem 9.}\\
\noi $\delta$ is a homomorphism in the {\em braided}\ tensor algebra
of operators.

For example:
\eq \delta(t_it_j)=\delta(t_i) \delta(t_j)
\en
indeed
\eq
\begin{array}{rcl}
\delta(t_i t_j)\dia a & \equiv & t_it_ja=t_it_j(a)+t_i(f_j{}^l*a)t_l\\
		      & = &t_i(t_j(a))
+(f_i{}^l*t_j(a))t_l+[t_i(f_j{}^l*a)]t_l+(f_i{}^m*f_j{}^l*a)t_mt_l\\
	       & = &  (t_i t_j \otimes 1 + (1 + \hat B)^{rs}_{ij}
		      (t_r \otimes t_s) + 1 \otimes t_i t_j)\dia a\\
	       & =  &(t_i \otimes 1 + 1 \otimes t_i)(t_j \otimes 1
		       + 1 \otimes t_j)\dia a\\
	       & =  &[\delta(t_i) \delta(t_j)]\dia a
\end{array}
\en
where we have used $\hat B_{ij}^{rs}\chi_r*f_s{}^l=f_i{}^l*\chi_j$
an identity easily derived from (\ref{dit}) and (\ref{fibif}),
or directly from $d(f_i{}^j*a)=f_i{}^j*da$.

As another example we have:
$$
\delta(\ell_V)=\delta i_V\,\delta d +\delta d\,\delta i_V
$$
indeed use the braiding (\ref{braid})
to show that
$(\delta i_V\,\delta d +\delta d\,\delta i_V)\dia\vart$
is given by (\ref{7.7}).\hfill \\

To summarize: we can state commutation relations for a general
differential operator $D$ and for any form $\vart$ in an abstract
way as
\eq
D \vart = \delta(D) \dia \vart,
\en
where $\delta$ is the {\em classical} coproduct on first order operators
and can be extended (Theorem 9) as a homomorphism to higher order
operators.

\section{A second Lie derivative and a puzzle}

It is natural to define a Lie derivative $\err_{h_i}$ of a generic
form $\tau$ along a
{\em right}-invariant vector field $h_i$
in terms of the {\em left} coaction $\Delta_\Gamma ~:$
$$
\err_h(\tau) \equiv (\chi \otimes i\!d)\Delta_\Gamma(\tau)
 = \tau * \chi~,
$$
just like it was natural that we used the {\em right} coaction,
when we defined $\ell_{t_i}$ in the previous section.
In this section we would like to compare the two definitions.

{}From the above definition we find
\eq
\begin{array}{rcl}
\err_{h_i}(\vart \vart') & = & \chi(\vart_1 {\vart'}_1) \vart_l {\vart'}_l\\
& = & \chi(\vart_1) \vart_l \vart' + f_i{}^j(\vart_1) \chi_j({\vart'}_1)
      \vart_l {\vart'}_l\\
& = & \err_{h_i}(\vart) \vart' + (\vart * f_i{}^j) \err_{h_j}(\vart')
\end{array}
\en
and, in a particular case:
\eq
\begin{array}{rcl}
\err_{h_i}(da\: b) & = & d(h_i(a)) b + (da * f_i{}^j) h_j(b)\\
    & = & d(h_i(a)) b + (da * f_i{}^j) S^{-1}N^k{}_j t_k(b).
\end{array} \label{errab}
\en
On the other hand we can use the results of the previous chapter to
give an alternative expression for the Lie derivative along the
right-invariant vector field $h_i$:
\eq
\begin{array}{lcl}
\ell_{h_i}& = &S^{-1}N^j{}_i \ell_{t_j} +
		 d(S^{-1}N^j{}_i)\wedge i_{t_j}\\
	  & = &M_i{}^j\ell_{t_j}+dM_i{}^j\wedge i_{t_j}\label{fine1}
\end{array}
\en
and hence
\eq
\begin{array}{rcl}
\ell_{h_i}(da\: b) & = & d(M_i{}^j t_j(a))b + M_i{}^j
    (f_j{}^k * da) t_k(b)\\
    & = & d(h_i(a)) b + M_i{}^j
    (f_j{}^k * da) t_k(b).
\end{array} \label{ellab}
\en
The difference between expressions (\ref{errab}) and (\ref{ellab}) is
a good index for the ``defect'' between left and right transports on
a quantum group:
\eq
\begin{array}{rcl}
(\ell_{h_i} - \err_{h_i})(da\: b)
& = & [M_i{}^j (f_j{}^k * da) - (da * f_i{}^j) M_j{}^k] t_k(b)\\
& = & -D\!I_i{}^k(a) t_k(b);
\end{array}
\en
\eq
D\!I_i{}^k(a) \equiv
[d(M_i{}^j)(f_j{}^k * a) - (a * f_i{}^j) d(M_j{}^k)]
\spz {\mbox{ (Defect Index)}}.
\en
In the last passage
we have used  the Leibniz rule for $d$ combined with (\ref{propM}).
The term in the square brackets is
always zero in the classical (undeformed) case. Note that
$(\ell_{h_i} - \err_{h_i})$
vanishes on $da$ and $b$ separately but not necessarily on $da\: b$.
The case of ``$b$'' confirms a result from
Section~\ref{secBBS}:
\eq
\ell_{h_i}(b) = S^{-1}N^j{}_i t_j(b) = h_i(b) = \err_{h_i}(b)
\en
and shows that we will not encounter any ambiguities or inconsistencies
as long as we deal with general vector fields and functions alone.
Problems can occur however when we start to introduce forms. Here
is what happens in a well known example:

In the approach \cite{Jurco,Zumino} to differential calculus on quantum groups
that is based on the $L^\pm$ and $T$ matrices of \cite{FRT} one can give the
following explicit expressions for the $f_i{}^j$ and the linearly independent
$\chi_k$ functionals
(using double-indices)
$$f_{i_1 i_2}{}^{j_1 j_2} = L^{+i_1}{}_{j_1} SL^{-j_2}{}_{i_2}
= \epsi\delta^{k_1}_{k_2}-L^{+k_1}{}_{j} SL^{-j}{}_{k_2}
$$
and for the adjoint representation $N^j{}_i,\:M_i{}^j$
$$N^{j_1 j_2}{}_{i_1 i_2} = S T^{i_1}{}_{j_1} T^{j_2}{}_{i_2}~~,~~~
 M_{i_1i_2}{}^{j_1j_2}=(S^{-1}T^{j_2}{}_{i_2})T^{i_1}{}_{j_1}~.$$
We want to investigate whether the
``Defect Index'' vanishes. First consider the $x^{a_1}{}_{a_2}$
of (\ref{chix})
such that $\chi^{k_1}{}_{k_2}(x^{a_1}{}_{a_2})
=\delta^{k_1a_1}\delta_{k_2a_2}$
and notice that $t^{k_1}{}_{k_2}(x_2^{a}{}_{a})S^{-1}(x_1^{a}{}_{a})=
\delta^{k_1}_{k_2}$
Therefore
\eq
(\ell_{h_i} - \err_{h_i})(da\:x_2^{a}{}_{a})S^{-1}(x_1^{a}{}_{a})
= \sum_{k_1=k_2} D\!I_i{}^k(a)\\
= \sum_{k_1=k_2}d(M_i{}^j)(f_j{}^k * a) - (a * f_i{}^j) d(M_j{}^k).
\en
Now $f_{j_1 j_2}{}^{k\: k}$
gives a non-trivial matrix
$Y^{j_1}{}_{j_2}$, while $M_{j_1 j_2}{}^{k\: k}$
will become the Kronecker symbol
$I\delta^{j_1}{}_{j_2}$ which vanishes as it is acted upon by $d$.
This is what remains if the Defect Index were zero:
$$
d(M_i{}^j) (Y_j * a) = 0
$$
for all $a \in A$ --- an incorrect
statement!

\section{Discussion}

We have seen that the bicovariant differential calculus studied so far
admits a natural contraction
operator along generic vector fields which leads to a Lie derivative
$\ell$ that is not unique.
We stress that there are no problems when we deal
with vector fields and functions only.
This should suffice in many physics applications (e.g. gauge theories)
that do not make direct reference to volume forms.
To evaluate the defect index we have used the Leibniz rule
for the exterior differential $d$ [see (\ref{fine1}) and (\ref{ellab})],
this shows that a possible way to obtain a differential calculus
with $\ell=\err$ is to define a new exterior differential $d$ that
satisfies a deformed Leibniz rule \cite{FP}.
The problem of the non-unique Lie derivative adds to another
mystery \cite{Su2} of differential calculi on quantum groups:
the number of forms is generally larger than one should expect from the
classical limit.
Notice also that the Defect Index is strictly connected with the
impossibility of the Lie derivative to be a good differential
operator i.e. satisfying a deformed Leibniz rule.
Indeed, consider $V_k\equiv M_k{}^jt_j$; since $V_k$ has a tensor index
we would like $\ell_{V_k}$ to exhibit a nice (braided)
deformation of Leibniz rule. Using 5) of Theorem 7 we find
\eq
\begin{array}{rcl}
\ell_{V_i}(a\vart')&=&(\ell_{V_i}a)\:\vart'+(a*f_i{}^j)\ell_{V_j}\vart'+
[dM_i{}^j\:f_j{}^k*a - a*f_i{}^jdM_j{}^k]\wedge i_{t_k}(\vart')\\
&=&(\ell_{V_i}a)\:\vart'+(a*f_i{}^j)\ell_{V_j}\vart'+
D\!I_i{}^k(a) \wedge i_{t_k}(\vart').
\end{array}
\en
It is the Defect Index that is responsible for the appearance of the
contraction operator in the Leibniz rule for the Lie derivative.
\sk
%

\vfill\eject

\begin{thebibliography}{99}

\bibitem{FRT} L.D. Faddeev,
N.Yu.\ Reshetikhin, L.A.\ Takhtajan,
{\sl Quantization of Lie Groups and Lie Algebras},
Algebra i Anal. 1 {\bf 1} (1989) 178 (Leningrad Math. J. {\bf 1} 193
(1990)).

\bibitem{Wor}
S.L.\ Woronowicz,
{\sl Differential Calculus on Compact Matrix Pseudogroups},
Commun. Math. Phys. {\bf 122} (1989) 125.

\bibitem{Jurco} B.\ Jur\v{c}o,
{\sl Differential calculus on quantized simple Lie groups},
Lett. Math. Phys. {\bf 22} (1991) 177.

\bibitem{Watamura}
U.\ Carow-Watamura, M.\ Schlieker, S.\ Watamura, W.\ Weich,
{\sl Bicovariant differential calculus on quantum groups
$SU_q(N)$ and $SO_q(N)$},
Commun. Math. Phys. {\bf 142} (1991) 605.

\bibitem{Castellani} L.\ Castellani,{\sl Differential calculus
on $ISO_q(N)$, quantum Poincar\'e algebra and $q$-gravity},
DFTT-70/93, hep-th 9312179, to be publ. in Commun. Math. Phys.;
{\sl The Lagrangian of q-Poincar\'e Gravity},
Phys. Lett. {\bf B327} (1994) 22.


\bibitem{SWZ1} P.\ Schupp, P.\ Watts and B.\ Zumino,
{\sl Differential Geometry on Linear Quantum Groups},
Lett. Math. Phys. {\bf 25} (1992) 139.

\bibitem{AC} P.\ Aschieri and L.\ Castellani,
{\sl An introduction to non-commutative differential geometry on quantum
groups},
Int. J. Mod. Phys. {\bf A8} (1993) 1667.

\bibitem{SWZ2} P.\ Schupp, P.\ Watts and B.\ Zumino,
{\sl Bicovariant Quantum Algebras and Quantum Lie Algebras},
 Commun. Math. Phys. {\bf 157} (1993) 305.

\bibitem{SWZ3} P.\ Schupp, P.\ Watts, B.\ Zumino,
{\sl Cartan Calculus on Quantum Lie Algebras},
XXII${}^{th}$ DGM, Ixtapa, 1993; Preprint LBL-34833 (1993);
hep-th/9312073.

\bibitem{SW} P.\ Schupp, P.\ Watts,
{\sl Universal and generalized Cartan calculus on Hopf algebras},
Preprint LBL-33655 (1994); hep-th/9402134.

\bibitem{Swe} M.E.\ Sweedler, {\sl Hopf Algebras}, (Benjamin,
New York, 1969).


\bibitem{Paolo} P.\ Aschieri,
{\sl The Space of Vector Fields on Quantum Groups},
preprint UCLA/93/TEP/25 (1993); hep-th/9311151.

\bibitem{Zumino} B.\ Zumino,
{\sl Introduction to the differential Geometry of Quantum Groups},
in: K. Schm\"udgen (ed.), Math. Phys. X, Proc. X-th IAMP Conf. Leipzig (1991),
Springer-Verlag (1992).

\bibitem{S1} P.\ Schupp,
{\sl Quantum Groups,
Non-Commutative Differential Geo\-me\-try and Appli\-ca\-tions},
PhD. the\-sis, Ber\-keley (1993); pre\-print LBL-34942;\\  hep-th/9312075.

\bibitem{S2} P.\ Schupp,
{\sl Cartan Calculus: Differential Geometry for Quantum Groups},
Enrico Fermi Summer School on Quantum Groups, Varenna (1994);
preprint LMU-TPW-94-8; hep-th/9408170.

\bibitem{Su} A.\ Sudbery,
{\sl Canonical Differential Calculus on Quantum General Linear Groups
and Supergroups},
Phys. Lett. {\bf B284} (1992) 61, erratum -ibid. {\bf B291} (1992) 519.

\bibitem{FP} L.D.\ Faddeev, P.N.\ Pyatov, {\sl The Differential Calculus
on Quantum Linear Groups}; hep-th/9402070

\bibitem{Su2} A.\ Sudbery, {\sl The Quantum Orthogonal Mystery}, 30${}^th$
Karpacz Winter School on Theoretical Physics: Quantum Groups, Poland (1994);
hep-th/9407110

\end{thebibliography}
\end{document}